\documentclass{amsart}         %Latex2e
\usepackage{amscd}             %Latex2e
\def\joinrel{\mathrel{\mkern-3mu}}
\def\relbar{\mathrel{\smash-}} 
\def\exlongarrow{\relbar\joinrel\relbar\joinrel\relbar\joinrel\rightarrow}

\newcommand{\Z}{{\mathbb Z}}
\newcommand{\C}{{\mathbb C}}
\newcommand{\R}{{\mathbb R}}
\newcommand{\Q}{{\mathbb Q}}
\renewcommand{\P}{{\mathbb P}}
\newcommand{\A}{\mathcal A}
\newcommand{\B}{\mathcal B}
\renewcommand{\H}{\mathcal H}

\newcommand{\bPt}{{\mathcal P}(\partial \H^3)}
\renewcommand{\O}{{\mathcal O}}

\renewcommand{\Im}{\operatorname{Im}}
\newcommand{\Gal}{\operatorname{Gal}}
\newcommand{\Ker}{\operatorname{Ker}}
\newcommand{\rank}{\operatorname{rank}}

\newtheorem{theorem}{Theorem}[section]
\newtheorem{thmA}{Theorem}
\newtheorem{thmB}{Theorem}
\newtheorem{thmBorel}{Borel's Theorem}

\newtheorem{proposition}[theorem]{Proposition}
\newtheorem{corollary}[theorem]{Corollary}

\newtheorem{conjecture}[theorem]{Conjecture}
\newtheorem*{conj}{Conjecture}
\newtheorem*{conjM}{Milnor's Conjecture}
\newtheorem*{conjR}{Ramakrishnan's Conjecture}

\theoremstyle{definition}

\newtheorem*{remark}{Remark}

\begin{document}

\title[Rationality problems for Chern-Simons invariants]{Rationality
Problems for
$K$-theory and Chern-Simons Invariants of Hyperbolic
3-Manifolds} 

\author{Walter D. Neumann}
\address{Department of Mathematics\\The University of
Melbourne\\Parkville, Vic 3052\\Australia}
\email{neumann@maths.mu.oz.au}
\author{Jun Yang}
\address{Department of Mathematics\\Duke University\\Durham NC 27707}
\email{yang@math.duke.edu}
\maketitle

%\noindent {\sl Mathematics Subject Classification:} 57M99; 19E99, 19F27.
%\bigskip
%\bigskip
 
This paper makes certain observations regarding some conjectures of
Milnor and Ramakrishnan in hyperbolic geometry and algebraic
$K$-theory. As a consequence of our observations, we obtain new
results and conjectures regarding the rationality and irrationality of
Chern-Simons invariants of hyperbolic 3-manifolds.

In this paper, by a {\em hyperbolic 3-manifold}, we shall mean a complete,
oriented hyperbolic 3-manifold with finite volume. So a
hyperbolic 3-manifold is compact or has finitely many cusps as its
ends (see e.g. \cite{thurston}). By a {\em cusped manifold} we shall mean
a non-compact hyperbolic 3-manifold. 

A hyperbolic 3-manifold $M$ is a quotient of the hyperbolic 3-space
$\H^3$ by a discrete subgroup $\Gamma$ of $PSL_2(\C)$ with finite
covolume. The isometry class of $M$ determines the discrete subgroup
up to conjugation.
To each subgroup $\Gamma$ of $PSL_2(\C)$, we can associate the
{\em trace field} of $\Gamma$, that is, the subfield of $\C$ generated
by traces of all elements in $\Gamma$. The trace field clearly depends
only on the conjugacy class of $\Gamma$, so one can define it to be
the trace field of the hyperbolic 3-manifold. However, the trace field
is not an invariant of commensurability class of $\Gamma$, although it
is not far removed. There is a notion of {\em invariant trace field}
due to Alan Reid (see \cite{reid}) which is a subfield of trace field
and does give an invariant of commensurability classes. Let
$\Gamma^{(2)}$ be the subgroup of $\Gamma$ generated by squares of
elements of $\Gamma$. The invariant trace field $k(M)$ of
$M$ is defined to be the trace field of $\Gamma^{(2)}$. We will
use invariant trace fields throughout this paper. We want to emphasize
here that each invariant trace field is a number field together with a
{\em specific embedding in $\C$}. 

Recall that a number field is totally real if all its embeddings into
$\C$ have image in $\R$ and totally imaginary if none of its
embeddings has image in $\R$. A {\em CM-field} is a number field with
complex multiplication, i.e., it is a totally imaginary quadratic
extension of a totally real number field.  A familiar class of
CM-fields is the class of cyclotomic fields. We will also use a
slightly more general notion. We say that an embedding $\sigma: F
\hookrightarrow \C$ of a number field $F$ is a {\em CM-embedding} if
$\sigma(F)$, as a subfield of $\C$, is an imaginary quadratic
extension of a totally real field. So $F$ is a CM-field if and only if
all its embeddings are CM-embeddings.

The Chern-Simons invariant of a compact $(4n-1)$-dimensional
Riemannian manifold is an obstruction to conformal immersion of the
Riemannian manifold in Euclidean space \cite{chern-simons}. For
hyperbolic 3-manifolds Meyerhoff \cite{meyerhoff} extended the
definition to allow manifolds with cusps. The Chern-Simons invariant
$CS(M)$ of a hyperbolic 3-manifold $M$ is an element in
$\R/\pi^2\Z$. It is {\em rational} (also called {\em torsion}) if it
lies in $\pi^2\Q/\pi^2\Z$. The main application of our paper to Chern
Simons invariants is the following theorem.
\begin{thmA}
The Chern-Simons invariant $CS(M)$
of a hyperbolic 3-manifold is rational if the associated embedding of
the invariant trace field $k(M)$ of $M$ is a CM-embedding.
\end{thmA}

Numerical evidence suggests that the Chern-Simons invariant of a
hyperbolic 3-manifold is usually irrational when the Theorem A
does not apply (though irrationality of Chern-Simons invariant --- and
of volume --- has not been proved for any example).  In particular,
our results lead to the following explicit irrationality conjecture.
\begin{conj} If the 
 invariant trace field $k=k(M)$ satisfies $k\cap\overline
k\subset\R$ then $CS(M)$ is irrational.  In particular, $CS(M)$ is
irrational if $k(M)$ has odd degree over $\Q$.
\end{conj}

The invariant trace field $k$ of an {\em arithmetic} hyperbolic
3-manifold has just one complex place (cf.~e.g.,
\cite{neumann-reid1}), so it either satisfies $k=\overline k$ and is
CM-embedded or it satisfies $k\cap\overline k\subset\R$. Thus, in this
case the conjecture would say that rationality or irrationality of the
Chern-Simons invariant is completely determined by whether or not
$k=\overline k$.  (Recently Reznikov \cite{reznikov} claimed
rationality of the Chern-Simons invariant for {\em any} compact
arithmetic hyperbolic 3-manifold.  However, his concept of
``arithmetic'' is what is usually called ``real-arithmetic'' and
implies that $k$ is CM-embedded, so his result follows from Theorem
A. But it also has an easier proof: a real-arithmetic manifold has a
cover with an orientation reversing isometry.)

Theorem A will follow from a more general result about Bloch
groups. In \cite{neumann-yang}, we showed that the Chern-Simons
invariant $CS(M)$ is determined modulo rationals by an element
$\beta(M)$ in the Bloch group $\B(k(M))$ (with $k(M)$ replaced by $\C$
this was known in the compact case by Dupont \cite{dupont}).

We will prove a general result on the Bloch group of a number field
$F$ which is embedded into $\C$ as a non-real subfield that is stable
under the complex conjugation.  If $F$ is a non-real Galois field, then any
embedding has this property.

Complex conjugation then induces an involution on $F$ and hence on
$\B(F)$.  Let $\B(F)\otimes\Q = \B_+(F) \oplus \B_-(F)$ be
the splitting into $\pm 1$ eigenspaces for such an involution. We
emphasize that even in the Galois case, the involution, and therefore
also this splitting, depends very much on the specific embedding of
$F$ into $\C$ that one uses (though if $F$ happens to be a CM-field,
then complex conjugation commutes with all the embeddings of $F$ into
$\C$ so the decomposition is independent of the choice of embedding).
The following theorem gives formulas for the ranks of $\B_+(F)$ and
$\B_-(F)$.  
\begin{thmB} Assume $F$ is a fixed non-real subfield of
$\C$ which is finite over $\Q$ and stable under complex conjugation in
$\C$. Let $r_2$ be the number of pairs of conjugate complex embeddings
of $F$, and let $r'_2$ be the number of pairs $\tau,\overline\tau$ of
complex embeddings of $F$ such that
$$
\tau(\overline x)=\overline\tau(x) \quad \forall x\in F\subset \C.
$$ 
Then the following equations hold:
\begin{gather*}
\rank(\B_-(F)) = \frac12(r_2+r'_2),\\ 
\rank(\B_+(F))=\frac12(r_2-r'_2).
\end{gather*}
\end{thmB}

In Theorem \ref{nonstable} we describe the situation when $F$ is not
stable under conjugation.\medskip

\thanks{We would like to thank Alan Reid for useful discussions. The
second author would also like to thank Professor John Millson for
useful discussion on Chern-Simons invariants and to thank Professor
Richard Hain and the referee for many useful comments and suggestions
on early drafts of the paper.}  \bigskip

\section{Hyperbolic 3-manifolds and Bloch groups}

We will use the upper half space model for the hyperbolic 3-space
$\H^3$. Note that $\H^3$ has a standard compactification by adding a
boundary $\partial \H^3=\C\cup \{\infty\} =
\P^1(\C)$. The group of orientation preserving isometries of $\H^3$ is
$PSL(2,\C)$, which acts on the boundary $\partial \H^3$ via linear
fractional transformations. Given any 4 ordered distinct points 
$(a_0,a_1,a_2,a_3)$ on $\partial \H^3$, there exists a unique element
in $PSL(2,\C)$ that maps $(a_0,a_1,a_2,a_3)$ to $(0, \infty, 1, z)$,
where  
$$z=\frac{(a_2-a_1)(a_3-a_0)}{(a_2-a_0)(a_3-a_1)} \in \C-\{0,1\},$$
is the cross ratio of $a_0, a_1, a_2, a_3$.
Given any {\em ideal} tetrahedron (vertices at the
boundary) in $\H^3$ and an ordering of its vertices, we can associate
to it the cross ratio of its vertices which is a complex number in
$\C-\{0,1\}$. If we change the ordering of the vertices by an even
permutation (so the orientation of the ideal tetrahedron is not
changed) we obtain 3 possible answers for the cross ratio: $z$,
$1-1/z$, and $ 1/(1-z)$.  

There are several definitions in the literature of the Bloch group
$\B(F)$ of an infinite field $F$. They are isomorphic with each other
modulo torsion of order dividing 6, and they are strictly isomorphic
for $F$ algebraically closed (see \cite{dupont-sah} and
\cite{neumann-yang} for a more detailed discussion).  The most widely
used are those of Dupont and Sah \cite{dupont-sah} and of Suslin
\cite{suslin}. We shall use a definition equivalent to Dupont and
Sah's, since it is the most appropriate for questions involving
$PGL_2(F)$. Suslin's definition is appropriate for questions involving
$GL_2(F)$ and gives a group which naturally embeds in ours with
quotient of exponent at most 2.  However, since torsion is irrelevant
to our current considerations, the precise definition we choose is, in
fact, unimportant.

Our $\B(F)$ is defined as follows. Let $\Z(\C-\{0,1\})$ denote the
free group generated by points in $\C-\{0,1\}$. Define the map
\begin{gather*}
\mu: \Z(\C-\{0,1\})\rightarrow \wedge^2_\Z(\C^\ast),\\
      \sum_i n_i [z_i] \mapsto \sum_i n_iz_i \wedge (1-z_i),
\end{gather*}
and let
$\A(F) = \Ker(2\mu)$.

The Bloch group $\B(F)$ is the quotient of $\A(F)$ by
all instances of the following ``5-term relation'' \marginpar{a
misprint in this equation in the published version was corrected here
03/96}
\begin{equation*}
 [x]-[y] + [y/x] - [\frac{1-x^{-1}}{1-y^{-1}}]+[\frac{1-x}{1-y}] =
0,
\end{equation*} 
where $x, y \in F -\{ 0, 1\}$. One checks easily that the map $\mu$
vanishes identically on this 5-term relation. 

\begin{remark} The above 5-term relation is that of 
Suslin \cite{suslin}, which differs slightly from the
following more widely used 5-term equation (first written down by
Bloch-Wigner) 
\begin{equation*}
 [x]-[y] + [y/x] - [\frac{1-y}{1-x}]+[\frac{1-y^{-1}}{1-x^{-1}}] =
0.
\end{equation*} 
These two different relations lead to mutually isomorphic $\B(\C)$ via
the map $[z]\mapsto[1/z]$. However, Suslin's 5-term relation seems
more suitable for general $F$ since it vanishes exactly under $\mu$
for an arbitrary $F$ while the 5-term equation of Bloch-Wigner only
vanishes up to 2-torsion\begin{footnote}{One could therefore give an
alternate definition of the Bloch group by replacing $2\mu$ by $\mu$
in the definition of $\A(F)$. This gives a group that lies between
Suslin's Bloch group and the one defined here. Although this may seem
more natural, the current definition is suggested by
homological considerations. If $F$ is algebraically closed the factor
$2$ makes no difference.}\end{footnote}.\end{remark}

The geometric meaning of the 5-term equation is the following:
\smallskip

Let $a_0,\ldots,a_4\in \P_1=\bPt$. Then the convex hull of these
points in $\H^3$ can be decomposed into ideal tetrahedra in two ways,
either as the union of the three ideal tetrahedra
$\{a_1,a_2,a_3,a_4\}$, $\{a_0,a_1,a_3,a_4\}$ and $\{a_0,a_1,a_2,a_3\}$
or as the union of two ideal tetrahedra $\{a_0,a_2,a_3,a_4\}$ and
$\{a_0,a_1,a_2,a_4\}$. Taking $(a_0,\ldots,a_4) = (0,\infty,1,x,y)$,
then the five terms in the above relation occur exactly as the five
cross ratios of the corresponding ideal tetrahedra.
\smallskip

We discuss briefly how a hyperbolic 3-manifold determines an element
in $\B(\C)$. For more details, see \cite{neumann-yang}.

It is straightforward for a closed hyperbolic 3-manifold $M$. Let $M =
\H^3/\Gamma$, where $\Gamma$ is a discrete subgroup of $PSL(2,\C)$.
The fundamental class of $M$ in $H_3(M,\Z)$ then gives rise to a class
in $H_3(PSL_2(\C);\Z)$ via the map
$$H_3(M, \Z)\cong H_3(\Gamma, \Z) \rightarrow H_3(PSL(2,\C),\Z).$$
Fix a base point in $\C-\{0,1\}$. Then via the natural action of
$PSL(2,\C)$ on $\P^1(\C)$, any 4-tuple of elements in $PSL(2,\C)$ 
gives a 4-tuple of elements in $\C$ which then gives rise to an
element in $\C$ via the cross ratio. This induces a well defined
homomorphism 
$$H_3(PSL(2,\C),\Z) \longrightarrow \B(\C).$$
Although the cross ratio depends on the particular choice of the base
point, the above homomorphism is independent of this choice (see
\cite[(4.10)]{dupont-sah}).  

In the cusped case, we note that by \cite{epstein-penner}, each cusped
hyperbolic 3-manifold $M$ has an ideal triangulation    
$$M = \Delta_1\cup\cdots\cup \Delta_n. $$
Let $z_i$ be the cross ratio of any ordering of the vertices of
$\Delta_i$ consistent with its orientation. The following relation
was first discovered by W. Thurston (unpublished --- see the remark in
Zagier \cite{zagier}; for a proof, see \cite{neumann-yang}):  
\begin{equation*}
\sum _i z_i \wedge (1-z_i) = 0,\quad\quad
\mbox{in}\;\wedge_{\Z}^2(\C^\ast). 
\end{equation*}
This relation is independent of the particular choice of the
ordering of the vertices of the $\Delta_i$'s. 

From the definition of Bloch groups, it is thus clear that an ideal
decomposition of a cusped hyperbolic 3-manifold determines an element
$\beta(M)=\sum_i[z_i]$ in the Bloch group $\B(\C)$. That such an element
is independent of the choice of the ideal triangulation is proved in
\cite{neumann-yang}.

In fact, what we proved in \cite{neumann-yang} is stronger. It
includes the following, which is what is needed for the purpose of
this paper.

\begin{theorem}\label{bloch-inv}
Every hyperbolic 3-manifold $M$ gives a
well defined element $\beta(M)$ in the Bloch group $\B(\C)\otimes \Q$
which is the image of a well defined element $\beta_k(M)$ in
$\B(k)\otimes\Q$ where $k$ is the invariant trace field of $M$.\qed
\end{theorem}

\section{Borel's Theorem}

The Bloch group of an infinite field $F$ is closely related to the
third algebraic $K$-group of $F$ (denoted by $K_3(F)$. There exists a
natural map
$$\B(F) \longrightarrow K_3(F)$$
due to Bloch \cite{bloch}. In particular, if $F$ is a number field,
then there is the isomorphism
$$\B(F)\otimes \Q \cong K_3(F)\otimes \Q.$$
For the precise relationship between the Bloch group and $K_3$ of an
infinite field, see the work of Suslin \cite{suslin}.

In \cite{borel}, Borel generalized the classical Dirichlet Unit
Theorem in classical number theory. The generalization started with
the observation that for a number field $F$, the unit group of the
ring of integers $\O_F$ is non-other than the first algebraic
$K$-group of $\O_F$, $K_1(\O_F)$. So the generalization is along the
line of higher $K$-groups. For the exact statement, we refer to
\cite{borel}. For the purpose of this paper, we state Borel's Theorem
for $K_3$ of a number field in terms of the Bloch group. The way we
state Borel's theorem here therefore incorporates a series of
important works by Bloch and Suslin.

Define the {\em Bloch-Wigner} function $D_2:
\C-\{0,1\}\rightarrow \R$ by (cf. \cite{bloch})  
$$D_2(z) = \Im \ln_2(z) + \log |z|\arg(1-z),\quad z\in \C -\{0,1\}$$
where $\ln_2(z)$ is the classical dilogarithm function.  The
hyperbolic volume of an ideal tetrahedron $\Delta$ with cross ratio
$z$ is equal to $D_2(z)$. It follows that $D_2$ satisfies the
five-term functional equation given by the 5-term relation, and
therefore $D_2$ induces a map
$$D_2: \B(\C) \longrightarrow \R,$$
by defining $D_2[z]=D_2(z)$.

Given a number field $F$, let $r_1$ and $r_2$ denote the number of real
embeddings $F\subset \R$ and the number of pairs of conjugate complex
embeddings $F\subset \C$ respectively. Let 
$\sigma_1,\ldots,\sigma_{r_1}, \sigma_{r_1+1},\bar{\sigma}_{r_1+1},
\ldots, \sigma_{r_1+r_2},\bar{\sigma}_{r_1+r_2}$ denote these
embeddings. Given a number field $F$, then one has a map
$$\begin{array}{rrcl}
c_2: & \B(F)& \longrightarrow & \R^{r_2}\\
     & \sum_i (n_i[z_i]) & \mapsto & (\sum_i n_i
D_2(\sigma_{r_1+1}(z_i)), \ldots, \sum_i n_i
D_2(\sigma_{r_1+r_2}(z_i))).
\end{array}
$$
\medskip

\begin{thmBorel} The kernel of $c_2$ is
exactly the torsion subgroup of $\B(F)$ and the image of $c_2$ is a
maximal lattice in $\R^{r_2}$. In particular, it follows that the rank
of $\B(F)$ is $r_2$.
\qed\end{thmBorel}

This theorem has some useful consequences.
Denote $\B(F)_\Q=\B(F)\otimes\Q$ for short. An immediate consequence
of the theorem is that an inclusion of number fields $F\hookrightarrow
E$ induces an injection $\B(F)/Torsion\to\B(E)/Torsion$ and hence an
injection $\B(F)_\Q\to \B(E)_\Q$.  Note that if $F\subset E$ is a
finite Galois extension of fields with Galois group $H$ then $H$ acts
on $\B(E)$.

\begin{proposition}\label{goodprops}
For any subfield $F$ of the number field $E$ identify $\B(F)_\Q$ with
its image in $\B(E)_\Q$.  Then if $F_1$ and $F_2$ are two subfields 
we have $$\B(F_1\cap F_2)_\Q=\B(F_1)_\Q\cap\B(F_2)_\Q$$
and if $E/F$ is a Galois extension with group $H$ then 
$$\B(F)_\Q=(\B(E)_\Q)^H$$
(the fixed subgroup of $\B(E)_\Q$ under $H$).
\end{proposition}

\begin{proof}
This result appears to be known to experts but we could find no
published proof, so we provide one.  It clearly suffices to prove the
results for the Bloch groups tensored with $\R$ rather than with $\Q$.
We first prove the Galois property in the case that $E$ is Galois over
$\Q$.  Let $G=\Gal(E/\Q)$.  If $\tau:E\to\C$ is our given embedding
and $\delta\in G$ the restriction of complex conjugation for this
embedding, then all complex embeddings have the form $\tau\circ
\gamma$ with $\gamma\in G$ and the conjugate embedding to $\tau\circ
\gamma$ is $\tau\circ\delta \gamma$.  Consider the map
$\B(E)\otimes\R\to\R{G}$ given on generators by
$$[z]\mapsto\sum_{\gamma\in G}D_2(\tau(\gamma(z)))\gamma.$$
By Borel's theorem this is injective with image exactly
$$\left\{\sum_{g\in G}r_{\gamma}\gamma \mid r_\gamma=-r_{\delta\gamma}
\text{ for all }\gamma\in G\right\}.$$
The $G$-action on $\B(E)$ corresponds to the action of $G$ from the
right on $\R G$.  Thus, if we identify $\B(E)\otimes\R$ with its image
in $\R G$, then $\B(E)^H\otimes\R$ is identified with the set of
elements $\sum r_\gamma\gamma\in\R G$ satisfying
$r_\gamma=-r_{\delta\gamma}$ and $r_\gamma=r_{\gamma\theta}$ for all
$\gamma\in G$ and $\theta\in H$.  That is, the function $r_\gamma$ is
constant on right cosets of $H$ and for the cosets $\gamma H$ and
$\delta\gamma H$ it is zero if they coincide and otherwise takes
opposite values on each. Thus the rank $r$ of $\B(E)^H\otimes\R$ is the
half the number of cosets $\gamma H$ for which
$\gamma^{-1}\delta\gamma\not\in H$.
Now, with $F=E^H$, the embedding $\tau\circ\gamma|F$ depends only on
the coset $\gamma H$ and is real or complex according as its
conjugation map
$\gamma^{-1}\delta\gamma$ does or does not lie in $H$, so the above
rank $r$ is
just $r_2(F)$.  Since the image of $\B(F)\otimes\R$ lies in
$\B(E)^H\otimes\R$ and has this rank, it must be all of
$\B(E)^H\otimes\R$. 

The case when $E$ is not Galois over $\Q$ now follows easily by
embedding $E$ in a larger field which is Galois over $\Q$.  Similarly,
for the intersection formula, by replacing $E$ by a larger field as
necessary we may assume $E$ is Galois over $F_1$ and $F_2$ with
groups $H_1$ and $H_2$ say. Then, identifying $\B(F_i)_\Q$ with its
image in $\B(E)_\Q$ we have $\B(F_1)_\Q\cap
\B(F_2)_\Q=\B(E)_\Q^{H_1}\cap\B(E)_\Q^{H_2} = \B(E)_\Q^{<H_1,H_2>}=
\B(E^{<H_1,H_2>})_\Q=\B(F_1\cap F_2)_\Q$.
\end{proof}

\begin{remark}
One can be a bit more precise about the functorial properties of $\B$.
In \cite{suslin} Suslin claims for his Bloch groups a transfer
homomorphism $L_{E/F}:\B(E)\to\B(F)$ with the usual properties of
transfer.  This implies that if the degree $[E:F]$ is $d$ then the map
$\B(F)\to\B(E)$ has kernel of exponent $d$ and its image lies in a
direct summand of $\B(E)$ with quotient of exponent $d$.  Since our
Bloch group may differ from Suslin's by exponent 2, we may have to
replace $d$ here by $2d$.
\end{remark}

\section{Proof and generalization of Theorem B}

We shall give a proof of Theorem B based directly on Borel's
theorem. A different proof can be given using Proposition
\ref{goodprops} and its consequence Theorem \ref{nonstable}.

We denote the complex conjugation in $\C$ by $\delta$. Let $F$ be a
fixed non-real subfield of $\C$ that is stable under complex
conjugation, i.e., $\delta(F) = F$. Assume $F$ is a finite extension
field of $\Q$. Let $r_2$ be as in Borel's Theorem. Then we may list
all the complex (non-real) embeddings of $F$ into $\C$ as $\tau_1,
\delta\tau_1, \ldots, \tau_{r_2}, \delta\tau_{r_2}$. Let $r'_2$ be the
number of conjugate pairs that commutes with $\delta$, i.e.,
$\tau_i\delta = \delta\tau_i$.  Renumbering if necessary, we may
assume $\tau_1,\ldots,\tau_l$ are the ones that commute with
$\delta$. Note that by our assumption on $F$, $r'_2$ is at least
one. The rest of the $\tau$'s won't commute with $\delta$, therefore
$\tau_i$ and $\tau_i\delta$, $i>r'_2$ will be in different conjugate
pairs (we use here that$F$ is non-real). So renumbering if necessary,
we may assume
$$
\tau_1,\ldots,\tau_{r'_2},
\tau_{r'_2+1},\tau_{r'_2+1}\delta,\ldots,\tau_m,\tau_m\delta
$$
gives exactly one representative from each conjugate pair of
embeddings of $F$ into $\C$, where $m=r'_2+(r_2-r'_2)/2$.

The complex conjugation on $F$ induces an involution on $\B(F)$. Let
$\B_+(F)$ and $\B_-(F)$ be the $\pm$ eigenspace of
$\B(F)_\Q=\B(F)\otimes\Q$. By Borel's Theorem, $\B(F)_\Q$ has a
$\Q$-basis $\alpha_1,\ldots,\alpha_{r_2}$. Let
$$u_i=\alpha_i-\delta\alpha_i,\quad v_i=\alpha_i+\delta\alpha_i,
1\leq i\leq r_2.$$
Then $u_i$'s and $v_i$'s span $\B_-(F)$ and $\B_+(F)$
respectively. Together, they span $\B(F)_\Q$. Hence by Borel's
Theorem, we know that the matrix 
$$
\left (
\begin{matrix}
D_2(\tau_1(u_1))&\cdots & D_2(\tau_m(u_1)) &
D_2(\tau_{r'_2+1}\delta(u_1)) & \cdots & D_2(\tau_{m}\delta(u_1)) \\ 
\vdots && \vdots & \vdots & & \vdots \\
D_2(\tau_{r_2}(u_{r_2}))& \cdots & D_2(\tau_m(u_{r_2})) & 
D_2(\tau_{r'_2+1}\delta(u_{r_2})) & \cdots &
D_2(\tau_{m}\delta(u_{r_2})) \\ 
D_2(\tau_1(v_1))&\cdots & D_2(\tau_m(v_1)) &
D_2(\tau_{r'_2+1}\delta(v_1))  & \cdots & D_2(\tau_{m}\delta(v_1)) \\ 
\vdots && \vdots & \vdots & & \vdots \\
D_2(\tau_{r_2}(v_{r_2}))& \cdots & D_2(\tau_m(v_{r_2})) & 
D_2(\tau_{r'_2+1}\delta(v_{r_2})) & \cdots &
D_2(\tau_{m}\delta(v_{r_2})) \\ 
\end{matrix}
\right )
$$
has rank $r_2$. Note that because the first $r'_2$ embeddings commute
with $\delta$, the entries of the last $r_2$ rows of the first $r_2$
columns are all 0's. Also, it follows from the equation
$\delta(u_i)=-\delta(u_i)$ and $\delta(v_i)=\delta(v_i)$, this matrix
has the following block form 
$$
\left (
\begin{matrix}
A_{r_2\times r'_2} & B_{r_2\times (r_2-r'_2)/2} & -B_{r_2\times (r_2-r'_2)/2} \\
0               & C_{r_2\times (r_2-r'_2)/2} &  C_{r_2\times (r_2-r'_2)/2} 
\end{matrix}
\right )
$$
So the matrix $A$ has to have rank $r'_2$. For the last $(r_2-r'_2)$-columns
to have rank $r_2-r'_2$, the matrices $B$ and $C$ must both have maximal
possible rank, that is, $(r_2-r'_2)/2$.  Since by Borel's Theorem, 
$$
\rank C = \rank \B_+(F),
$$
and $\rank \B_+(F) + \rank \B_-(F) = r_2$, Theorem B follows.\qed\smallskip

We can also describe the situation when $F\subset \C$ is a number
field that is not stable under conjugation.  If $E\subset\C$ is any
number field containing $F$ with $E=\overline E$ then $\B_+(E)$ and
$\B_-(E)$ are defined, so we can form
$$\B_+(F):=\B_+(E)\cap\B(F)_\Q\quad\text{and}\quad 
 \B_-(F):=\B_-(E)\cap\B(F)_\Q.$$
These subgroups are independent of the choice of $E$, but in general
they will not sum to $\B(F)_\Q$.

Denote $F_\R=F\cap\R$ and let $F'=F\cap\overline F$.
Clearly $F'$ contains $F_\R$, and $F_\R$
must be the fixed field of conjugation on $F'$.  Thus either
$F'=F_\R$ or $F'$ is an imaginary quadratic extension of $F_\R$.
Now $F'$ is a field to which Theorem B applies, so
$\B(F')_\Q=\B_+(F')\oplus\B_-(F')$, with the ranks of the summands
given by Theorem B.

\begin{theorem}\label{nonstable}
$\B_-(F)=\B_-(F')$ and $\B_+(F)=\B_+(F')=\B(F_\R)_\Q$.
\end{theorem}

\begin{corollary}\label{corNS}
$\B_+(F)$ is trivial if and only if $F_\R$ is totally real.

$\B_-(F)$ is trivial if and only if $F'=F_\R$.

$\B_-(F)=\B(F)_\Q$ if and only if $F=F'$ and $F_\R$ is totally real;
that is either $F$ is totally real or the embedding
$F\hookrightarrow\C$ is a CM-embedding.
\end{corollary}

\begin{proof}[Proof of Theorem \ref{nonstable}]
We work in a Galois superfield $E$ of $F$ and identify Bloch groups
with their images in $\B(E)_Q$.  Let $G=\Gal(E/\Q)$, so
$H=\Gal(E/F)\subset G$ is the subgroup which fixes $F$.  We fix an
embedding $E\subset\C$ extending the given embedding of $F$ and denote
complex conjugation for this embedding by $\delta$. Then the subgroup
$H_\R$ generated by $H$ and $\delta$ is $\Gal(E/F_\R)$, so it follows
from Proposition \ref{goodprops} that
$\B_+(F)=\B(E^H)_\Q\cap\B(E)_\Q^\delta=
(\B(E)_\Q^H)^\delta=\B(E)_\Q^{H_\R}=\B(F_\R)_\Q$.  Moreover, $\B_-(F)$
is fixed by both $H$ and $\delta H\delta$ and hence by the group $H'$
that they generate.  But $H'$ is the Galois group in $E$ of
$F\cap\delta(F)=F\cap\overline F=F'$.  Thus $\B_-(F)$ is in
$\B(E)_\Q^{H'}= \B(F')_\Q$.  We thus obtain inclusions
$\B_{\pm}(F)\subset\B_\pm(F')$, and the reverse inclusions are trivial.
\end{proof}

\begin{proof}[Proof of Corollary]
The first two statements follow immediately from Theorem
\ref{nonstable} and Theorem B.  For the third, note that
$\B_-(F)\subset\B(F')_\Q$ and if $F'\ne F$ then $F$ has strictly more
complex embeddings than $F'$ so $\B(F')_\Q\ne\B(F)_\Q$.  Thus, to have
$\B_-(F)=\B(F)_\Q$ we must have $F=F'$.  The claim then follows
directly from Theorem B.
\end{proof}

\begin{remark}
We have pointed out at the beginning of sect.~2 that $\B(F)$ could
have been replaced by $K_3(F)$ in all our discussions.  The analog of
Borel's theorem holds for $K_i(F)$ for all $i\equiv3$ (mod $4$), so
the results described above are also valid for these $K$-groups.  When
$1<i\equiv 1$ (mod $4$) Borel's theorem gives a map
$K_i(F)\to\R^{r_1+r_2}$ whose kernel is torsion and whose image is a
lattice. The only change is that one obtains
$r_1+\frac12(r_2+r_2')$ and $\frac12(r_2-r'_2)$ as
the dimensions of the $+$ and $-$ eigenspaces in the analog of Theorem
B, and Corollary \ref{corNS} therefore also needs modification.  We
leave the details to the reader. The basic point is that if
$E\subset\C$ is Galois over $\Q$ with group $G$ and $\delta$ is its
conjugation then $K_i(E)\otimes\R$ is $G$-equivariantly isomorphic to
$\{\sum r_\gamma\gamma\in\R G |
r_\gamma=(-1)^{(i-1)/2}r_{\delta\gamma}\}$ for $i>1$ and odd.
\end{remark}

\section{Milnor's and Ramakrishnan's conjectures}

Milnor \cite{milnor1} made the following conjecture motivated by the
fact that $D_2(z)$ represents the volume of an ideal tetrahedron.  For
the significance of this conjecture in hyperbolic geometry and number
theory, see \cite{milnor1}, \cite{milnor2}.

\begin{conjM}
For each integer $N\geq 3$, the real
numbers $D_2(e^{2\pi \sqrt{-1}j /N})$, with $j$ relatively prime to $N$
and $0 < j < N/2$, are linearly independent over the rationals.
\end{conjM}

A field homomorphism $\tau: F\rightarrow K$ clearly induces a
homomorphism on the Bloch groups $\B(F)\rightarrow \B(K)$ which, by
abuse of notation, will again be denoted by $\tau$. 

Given a cyclotomic field $F=\Q(e^{2\pi \sqrt{-1}/N})$, the elements
$[e^{2\pi i j/N}]$, with $j$ relatively prime to $N$ and $0<j<N/2$,
form a basis of the Bloch group $\B(F)\otimes \Q$ (see Bloch
\cite{bloch}). Hence Milnor's conjecture can be reformulated that
$D_2:\B(F)\to \R$ given on generators by $[z]\mapsto D_2(z)$ is
injective for a cyclotomic field $F$.

Note that for a general number field $F$ the above map $D_2$
vanishes on $\B_+(F)$.  By Corollary \ref{corNS}, the following is
thus the strongest generalization of Milnor's conjecture that one
might hope for.

\begin{conjecture}\label{conjB}
If $F\subset\C$ is a number field with $F\cap\R$ totally real then the
map $D_2:\B(F)\to \R$ given by $[z]\mapsto D_2(z)$ is injective.
\end{conjecture}

In the special case that $F$ is Galois over $\Q$ the condition
$F\cap\R$ totally real says that $F$ is a CM-field. In this case we
have the following proposition (cf. Prop. 7.2.5 of
\cite{ramakrishnan})

\begin{proposition}
Suppose that $F$ is a Galois CM-field over $\Q$. If for one complex
embedding $\tau: F\hookrightarrow \C$, the map $D_2\circ\tau$
is injective (that is, Conjecture \ref{conjB} holds), then it is for
all complex embeddings.
\end{proposition}

\begin{proof}
Let $\rho: F\hookrightarrow \C$ be another embedding of $F$. There
exists $\gamma\in \Gal(F/\Q)$ such that $\rho = \tau\gamma$. Let
$\omega\in B(F)\otimes \Q$.  If
$$D_2\circ
\rho(\omega)= D_2\circ\tau(\gamma(\omega)) = 0,$$
 then by the injectivity of $D_2\circ\tau$, $\gamma(\omega)=0$.
Since $\gamma: \B(F)\rightarrow \B(F)$ is clearly an isomorphism, it
follows that $\omega = 0$ in $\B(F)\otimes \Q$. Hence, $D_2\circ\rho$
is also injective.
\end{proof}

The map $D_2: \B(\C)\rightarrow \R$ is the imaginary part of the more
general {\em Bloch} map
$$\rho: \B(\C) \longrightarrow \C/\Q(2),$$
where the notation $\Q(k)$ denotes the subgroup $(2\pi \sqrt{-1})^k\Q$
of $\C$. The definition of $\rho$ is given as follows:

For $z\in \C-\{0,1\}$, define
$$\rho(z) = \log z\wedge\log (1-z) + 2\pi \sqrt{-1}\wedge
\frac{1}{2\pi \sqrt{-1}}(\ln_2(1-z)-\ln_2(z) - \pi^2/6) \in \Lambda^2_\Z \C.$$
See section 4 of \cite{dupont-sah} or \cite{hain} for the meaning of
this map and for further details. The exact formulae given here is
borrowed from \cite{hain}.  This map obviously induces a map
$$\rho: \A(\C)\longrightarrow \wedge^2_\Z\C.$$
It turns out that $\rho$ vanishes on the 5-term relation, hence it
induces a map
$$\rho: \B(\C)\longrightarrow \wedge^2_\Z\C.$$
Finally, it follows from the fact that every element $\alpha$ of
$\B(\C)$ satisfies the relation $\mu(\alpha)=0$ that the image of this
last map lies in the kernel of the map
$$e: \wedge^2_Z \C
\stackrel{\wedge^2\exp}{\exlongarrow} \wedge^2_\Z\C^\ast.$$
The kernel of $e$ is $\C/\Q(2)$. Hence this induces the Bloch map.

Ramakrishnan \cite{ramakrishnan} generalized Milnor's conjecture in
the following form\footnote{Both Ramakrishnan's conjecture and
Milnor's original conjecture are more general in that they apply to
all the odd-degree higher $K$-groups. We refer to \cite{ramakrishnan}
for more details. Likewise, Conjecture \ref{conjB} can be stated for
higher $K$-groups in similar fashion.}.

\begin{conjR}
 For every subfield $F\stackrel{\sigma}{\hookrightarrow} \C$,
the map  
$$\B(F)\otimes\Q\stackrel{\sigma}{\longrightarrow}\B(\C)\otimes \Q 
\longrightarrow \C/\Q(2)$$
is injective. 
\end{conjR}

The Bloch-Wigner function $D_2$ is the imaginary part of the Bloch map
$\rho$, and it vanishes identically on $\B_+(k)$. On the other hand,
it follows from a routine calculation that the real part of the Bloch
map vanishes identically on $\B_-(k)$.  

In particular, $\rho$ just
reduces to $D_2$ if $\B_-(F)=\B(F)_\Q$.  By Corollary \ref{corNS} we
thus have

\begin{proposition}
If $F\subset\C$ is a CM-embedded field, then the
Ramakrishnan Conjecture for this particular embedding of $F$ is
equivalent to Conjecture \ref{conjB}.\qed
\end{proposition}

On the other hand 
\begin{proposition}
The truth of the Ramakrishnan Conjecture for a field $E=\overline
E\subset \C$ would imply Conjecture \ref{conjB} for any subfield of
$E$.
\end{proposition}

\begin{proof}
Since the real and imaginary parts of $\rho$ vanish on $\B_-(E)$ and
$\B_+(E)$ respectively, the Ramakrishnan conjecture for $E$ is
equivalent to the conjecture that the kernel of the real part of
$\rho$ is exactly $\B_-(E)$ and the kernel of the imaginary part, that
is $Ker(D_2)$, is exactly $\B_+(E)$. Thus $D_2$ would have zero kernel
on $\B(F)_\Q$ for any subfield $F$ of $E$ satisfying
$\B(F)_\Q\cap\B_+(E)=\{0\}$, which is equivalent to the condition of
Conjecture \ref{conjB} by Corollary \ref{corNS}.
\end{proof}

\section{Chern-Simons invariants and the regulator map} As we have
seen from the above discussion, hyperbolic 3-manifolds and their
volumes all have interesting $K$-theoretic interpretation. Another
important invariant in the theory of hyperbolic 3-manifolds is the
Chern-Simons invariant of \cite{chern-simons} and
\cite{meyerhoff}. What is its relation with $K$-theory?  As discussed
in sect.~1, given hyperbolic 3-manifold $M$ represents a class
$\beta(M)$ in $\B(\C)$. For $M$ compact, the next theorem follows from
Theorem 1.11 of Dupont \cite{dupont}. A proof for cusped $M$ can be
found in \cite{neumann-yang}.

\begin{theorem}\label{formula}
The Chern-Simons invariant of 
$M$ mod $\Q(2)$ is equal to the real part of the image of $\beta(M)$
under the regulator map   
$$\rho: \B(\C)\longrightarrow \C/\Q(2).\eqno{\qedsymbol}$$
\end{theorem}

Now given a hyperbolic 3-manifold $M$ with invariant trace field $k$
and the associated embedding $\sigma: k\hookrightarrow \C$, the class
$\beta(M)\in \B(\C)$ associated to $M$ is in the image of
$\B(k)_\Q\stackrel{\sigma}{\rightarrow} \B(\C)$ (Theorem
\ref{bloch-inv}). If $\sigma$ is a $CM$-embedding, then it follows
from Corollary \ref{corNS} that the real part of $\rho$ is trivial
on $\B(F)\otimes \Q$. Theorem A in the introduction now follows
immediately from Theorem \ref{formula}. Similarly, as in the previous
section, Corollary \ref{corNS} implies that the irrationality
conjecture for Chern Simons invariant of the Introduction would be
implied by the Ramakrishnan Conjecture.
 
We have seen that the imaginary part of $\rho$ is essentially the
volume map while the real part of $\rho$ can be called the
Chern-Simons map. Reinterpreting the discussion of the previous
section, this means that if $k=\overline k\subset \C$, then the volume
map on $\B(k)$ factors through $\B_-(k)$ and Chern-Simons map on
$\B(K)$ factors over $\B_+(k)$.  By Theorem B we thus get bounds of
$\frac12(r_2+r'_2)$ and $\frac12(r_2-r'_2)$ on the number of
rationally independent volumes resp.\ Chern-Simons invariants for
manifolds having invariant trace field contained in our given
$k$. Ramakrishnan's Conjecture says the image of $\rho$ has $\Q$-rank
$r_2$.  This is equivalent to the conjecture that the $\Q$-ranks of
the images of $vol\colon \B(k)\to \R$ and $CS\colon\B(k)\to \R/\Q$ are
exactly $\frac12(r_2+r'_2)$ and $\frac12(r_2-r'_2)$.

\bibliographystyle{plain}

\begin{thebibliography}{99}

\bibitem{beilinson}
A.~Beilinson: {\it Higher regulators and values of $L$-functions,}
(English translations), Journal of Soviet Math. 30, no.~ 2 (1985),
2036--2070. 

\bibitem{bloch}
S.~Bloch: {\it Higher regulators, algebraic $K$-theory, and zeta
functions of elliptic curves,} Lecture notes U.C. Irvine (1978).

\bibitem{borel}
A.~Borel: {\it Cohomologie de $SL_n$ et valeurs de fonction zeta aux  
points entiers}, Ann. Sci. Ecole Norm. Sup. (4) 7, (1974), 613--636.

\bibitem{chern-simons}
S.~Chern, J.~Simons, {\it Some cohomology classes in principal fiber
bundles and their application to Riemannian geometry},
Proc. Nat. Acad. Sci. U.S.A. 68 (1971), 791-794.

\bibitem{dupont}
J.~Dupont: {\it The dilogarithm as a characteristic class for flat  
bundles},
J.\ Pure and App.\ Algebra 44 (1987), 137--164.

%\bibitem{dupont2}
%J.~Dupont: {\it Algebra of polytopes and homology of flag complexes},
%Osaka J.\ Math.\ 19 (1982), 599--641.

\bibitem{dupont-sah}
J.~Dupont, H.~Sah: {\it Scissors congruences II}, J.\ Pure and App.\
Algebra 25 (1982), 159--195.

\bibitem{epstein-penner}
D.~B.~A.~Epstein, R.~Penner, {\it Euclidean decompositions of
non-compact hyperbolic manifolds,} J.\ Diff.\ Geom. 27 (1988), 67-80.

%\bibitem{gillet}
%H.~Gillet: {\it Riemann-Roch theorems for higher algebraic
%$K$-theory,} Adv.\ Math.\ 40, no. 3 (1981), 203--289.

\bibitem{hain}
R.~Hain: {\it Classical polylogarithms}, Motives, Proc. Symp. Pure
Math, to appear.

%\bibitem{hain-mac}
%R.~Hain, R.~MacPherson: {\it Higher Logarithms}, Illinois J.\ Math.\  
%34 (1990), 392--475.
%
%\bibitem{macbeath}
%A.~M.~Macbeath: {\it Commensurability of cocompact three-dimensional
%hyperbolic groups.} Duke Math. J. 50(1983), 1245-1253.
%
\bibitem{meyerhoff}
R.~Meyerhoff: {\it Hyperbolic 3-manifolds with equal volumes but
different Chern-Simons invariants}, in Low-dimensional topology and
Kleinian groups, edited by D.\ B.\ A.\ Epstein, London Math.\ Soc.\
lecture notes series, 112 (1986) 209-215.

\bibitem{milnor1}
J.~Milnor: {\it Hyperbolic geometry: the first 150 years}, Bulletin
AMS 6 (1982), 9-24.
 
\bibitem{milnor2}
J.~Milnor: {\it On polylogarithms, Hurwitz zeta functions, and their
Kubert identities,} L'Einseignement Math., t. 29 (1983), 281-322.

\bibitem{neumann}
W.~D.~Neumann: {\it Combinatorics of triangulations and the Chern
Simons invariant for hyperbolic 3-manifolds},
in {\it Topology 90, Proceedings of the
Research Semester in Low Dimensional Topology at Ohio State} (Walter
de Gruyter Verlag, Berlin - New York 1992), 243--272.

\bibitem{neumann-reid1}
W.~D.~Neumann, A.~W.~Reid: {\it Arithmetic of hyperbolic manifolds},
{\it Topology 90, Proceedings of the Research Semester in Low
Dimensional Topology at Ohio State} (Walter
de Gruyter Verlag, Berlin - New York 1992), 273--310.

\bibitem{neumann-yang}
W.~D.~Neumann, J.~Yang: {\it Bloch invariants of hyperbolic
3-manifolds}, in preparation.  

\bibitem{neumann-zagier}
W.~D.~Neumann, D.~Zagier: {\it Volumes of hyperbolic 3-manifolds},
Topology 24 (1985), 307-332.
 
\bibitem{ramakrishnan}
D.~Ramakrishnan: {\it Regulators, algebraic cycles, and values of 
$L$-functions}, Contemp.\ Math.\ 83 (1989), 183--310.

%\bibitem{rapoport}
%M.~Rapoport: {\it Comparison of the regulators of Beilinson and of
%Borel,} in ``Beilinson's Conjectures on Special Values of
%$L$-functions'', edited by M. Rapoport et al. Perspectives in
%Mathematics. Vol 4, Academic Press, 1988, 169--192.

\bibitem{reid}
A.~W.~Reid: {\it A note on trace-fields of Kleinian groups.} Bull.
London Math. Soc. 22(1990), 349-352.

\bibitem{reznikov}
A.~G.~Reznikov: {\it Rationality of secondary classes}, Preprint 26,
The Hebrew University (1993).

%\bibitem{sah}
%C.~S.~ Sah: {\it Scissors congruences, I, Gauss-Bonnet map}, Math.\
%Scand.\ 49 (1982) 181-210.
 
\bibitem{suslin}
A.~A.~Suslin: {\it Algebraic $K$-theory of fields,} Proc.\ Int.\ Cong.\
Math.\ Berkeley 86, vol. 1 (1987), 222-244.

\bibitem{thurston}
W.~P.~Thurston: {\it The geometry and topology of 3-manifolds.}
Lecture notes, Princeton University, 1977.

\bibitem{zagier}
D.~Zagier: {\it The Bloch-Wigner-Ramakrishnan polylogarithm function,}
Math. Annalen 286 (1990), 613--624.

\end{thebibliography}

\end{document}